# A Machine Learning Model for Solving Lane-Emden Equation using Legendre Wavelet Neural Network

Vijay Kumar Patel[1], Abhishekh[2], Dileep Kumar[3*], and Nitin Kumar[4]

[1] Department of Mathematics and Computational Sciences, Indian Institute of Information Technology Surat, India

[2,3] School of Advanced Sciences and Languages, VIT Bhopal University, Bhopal-Indore Highway, Kothri-Kalan, Sehore, MP, India

[4] A.S.K.S. Autonomous P.G. College, Allahabad State University, Fatehpur-212656, India
vijaybhuiit@gmail.com, abhibhu1989@gmail.com, dileepkumar0204@gmail.com, nitin-shanu1@gmail.com

*Corresponding Author: dileepkumar0204@gmail,com

**Abstract.** As we know differential equations are very useful for electrical engineers to solve a variety of problems like: voltage across a capacitor, input versus output voltage, etc. Therefore, the goal of this paper is to find the solutions of non-linear differential equations based on the Lane–Emden equation of second order using the Legendre wavelet neural network (LWNN) method. Here all the considered equations are singular initial value problems. To manage the singularity challenge, we have employed an artificial neural network method. This approach utilizes a neural network of a single layer, where the hidden layer is omitted by enlarging the input using Legendre wavelets functions. We have applied a feed-forward neural network method to the proposed problem along with the principle of error backpropagation. The effectiveness of the Legendre wavelet Neural Network method is validated through Lane–Emden equations..

## 1. Introduction

We know that differential equations form the foundation of physical systems. Numerous models in mathematics, engineering, economics, and physics can be represented by ordinary differential equations or partial differential equations. In many instances, obtaining exact solutions to these equations is challenging. And finding exact solutions to nonlinear differential equations can also be challenging or even impossible in many cases. As a result, various analytical, numerical, and approximation methods are employed to solve them, such as perturbation techniques, finite difference methods, finite element methods, variational methods, Runge–Kutta, and predictor-corrector, operational matrix method [1-13], have been stablished to solve them. All the numerical approaches involve discretizing the domain into a finite number of grids/messes or regions where all the functions are locally approximated.
Now a days, Artificial Neural Networks (ANN) [14] have been applied to find the solution of the initial and boundary value problems using various machine learning techniques. ANN-based numerical solutions offer several opportunities. The trial solutions generated by ANN use only a single independent variable, in spite of the problem's dimensionality. These numerical solutions are continuous throughout the entire domain. Additionally, unlike few other methods that are typically iterative and require a fixed step size before computation, ANN does not require such constraints. In traditional methods, if we need the solution between steps after obtaining the initial result, the process must be repeated from the beginning. ANN provides a solution to this issue by allowing us to obtain results at arbitrary point without repeating the entire iterative process. Furthermore, it can be used to produce numerical results directly using a black box. Therefore, this manuscript consider a the following general form of Lane–Emden equations which will solving using a shifted Legendre Neural Network based model for solving.

The general form of Lane-Embden equation [15] is

$$\frac{d^2\xi}{d\eta^2} + \frac{2}{x}\frac{d\xi}{d\eta} + g(\eta, \xi) = f(\eta), \eta \geq 0, \tag{1}$$

with $\quad \xi(0) = \xi_0, \xi'(0) = 0,$
where $g(\eta, \xi)$ is a nonlinear function of $\eta$ and $\xi$ .

The equation (1) has many applications in engineering and physics. Most of the applications are depend on $g(\eta, \xi)$ .

## 2. Preliminaries

In this part, we have illustrated structure of single layer Legendre neural network model and also we have described its learning algorithm.

### 2.1. Legendre wavelets

In this subsection, we utilize Legendre wavelet neural network (LWNN) based on a single-layer model for this problem. The LWNN structure comprises node of a single input, an expandable functional block based on Legendre polynomial functions, and a node of single output. The neural model's architecture is decompose into two parts: the computational transformation part and the learning part. In the computational transformation part, all input datum is developed into multiple basis terms using the Legendre polynomial functions. It is efficiently creating a new input vector. Let's denote the input data as $\eta = (\eta_1, \eta_2, \ldots, x_h)^T$, where $\eta$ is the single input node contains $h$ data points. The Legendre wavelets functions are serve as the orthogonal basis for this expansion. The Legendre wavelets defined as follows: The first two Legendre polynomials on the interval $[-1,1]$ are well-known as $P_0(x) = 1$, $P_1(x) = x$.

The general order Legendre polynomial functions on the interval $[-1,1]$ can be find by the following recursive formula [1]

$$(k+1)P_k(x) = (2k+1)x\,P_k(x) - k\,P_{k-1}(x), \tag{2}$$

where $P_k(x)$ denotes nth order Legendre polynomial.

The Legendre wavelets [4] defined as follows:

$$L_{lm}(x) = \begin{cases} \sqrt{m+\frac{1}{2}}\,2^{k/2} P_m(2^k - 2l + 1), & \text{if } \frac{l-1}{2^{k-1}} \leq x < \frac{l}{2^{k-1}}, \\ 0, & \text{otherwise} \end{cases} \tag{3}$$

where $k$ is natural numbers, $l = 1,2,3\ldots,2^{k-1}$ and $m = 0,1,2,\ldots,M-1$.

In this approach, an $n$-dimensional input is transformed into an $m$-dimensional set of Legendre wavelets function. The benefit of using the Legendre wavelet Neural Network (LWNN) lies in obtaining results with network of a single-layer, achieved by increasing the dimension of input through Legendre wavelets function. The network's architecture, featuring the first five Legendre wavelet function along with input and output layer with a single node, is calculated by using the value of $k$

### 2.1. An algorithm of machine learning based on LWNN

An algorithm is employed to update the parameters of network and minimize the error. In this case, the algorithm of back propagation for error is utilized to adjust the weights of the LWNN. This involves calculating error for the gradient with parameters $p$. The function $tanh(\eta)$ is used for the activation function. The output of the network with the given input $\eta$ and $p$, can then be calculated as follows:

$$N(\eta, p) = tanh(C).$$

Here $C$ is the linear sum can be written as

$$C = \sum_{j=1}^{n} \sum_{i=1}^{m} w_{ji} L_{j-1,i-1}(\eta),$$

where $\eta = (\eta_1, \eta_2, \ldots, \eta_h)^T$ is the input, $w_i$ and $L_{j-1,i-1}(\eta)$ denote the expanded weight vectors and input date respectively of the LWNN.

Now, the weight of LWNN can be defined as

$$w_{ji}^{k+1} = w_{ji}^{k} + \Delta w_{ji}^{k} = w_{ji}^{k} - \rho \frac{\partial E(\eta, p)^k}{\partial w_{ji}^{k}}, \tag{4}$$

where the iteration $k$ is used to revise the weights, $\rho$ is learning parameter and $E(\eta, p)$ denotes the error function.

## 3. Neural network method based on Legendre wavelets

In this part, we will present a Legendre wavelet neural network (LWNN) method. Let $\xi_t(\eta, p)$ be the trial solution of LWNN with $p$ then it can be written as

$$\xi_t(\eta, p) = G(\eta) + H(\eta, N(\eta, p)), \tag{5}$$

where $G(\eta)$ is the first term, which is not included and adjusted any parameters and is designed solely to satisfy the boundary or initial conditions. In contrast, the second term, $H(\eta, N(\eta, p))$ is output of the LWNN with the single network, which depends on adjustable parameters $p$ and the input $\eta$ . As previously noted, the LWNN model consists of a single layer with single input node and single output node $N(\eta, p)$, where $\eta$ is the input and $p$ represents the network parameters and

$$N(\eta, p) = tanh(C), \quad \text{where } C = \sum_{j=1}^{n} \sum_{i=1}^{m} w_{ji} L_{j-1,i-1}(\eta).$$

Now, let us focus on the expression specifically for ordinary differential equations of second-order, as our goal is to solve the second order Lane–Emden equation. The proposed problem under consideration is expressed as follows:

$$\frac{d^2\xi}{d\eta^2} = f(\eta, \xi, \frac{d\xi}{d\eta}), \eta \in [a, b], \tag{6}$$

with the initial conditions $\xi(a) = G, \xi'(a) = G'$.

The LWNN trail solution is

$$\xi_t(\eta, p) = G + G'(\eta - a) + (\eta - a)^2 N(\eta, p). \tag{7}$$

Here $\eta$ is input and $N(\eta, p)$ is the output, $p$ is parameter and $\xi_t(\eta, p)$ is the trial solution which are satisfying the initial conditions.

Now, $N(\eta, p)$ needs to be minimized can be expressed as follows:

$$E(\eta, p) = \sum_{i=1}^{h} \frac{1}{2}\left(\frac{d^2\xi_t(\eta_i, p)}{d\eta^2} - f\left[\eta_i, \xi_t(\eta_i, p), \frac{d\xi_t(\eta_i, p)}{d\eta}\right]\right)^2. \tag{8}$$

The weights between the input and output layers are adjusted according to the following rule:

$$w_j^{k+1} = w_j^k + \Delta w_j^k = w_i^k - \rho \frac{\partial E(\eta, p)^k}{\partial w_j^k}. \tag{9}$$

Here

$$\frac{\partial E(\eta, p)}{\partial w_j} = \frac{\partial}{\partial w_j}\left(\sum_{i=1}^{h} \frac{1}{2}\left(\frac{d^2\xi_t(\eta_i, p)}{d\eta^2} - f\left[\eta_i, \xi_t(\eta_i, p), \frac{d\xi_t(\eta_i, p)}{d\eta}\right]\right)^2\right). \tag{10}$$

Finally, the converged results from the LWNN can be substituted into Eq. (7) to obtain the solutions.

## 4. Numerical results and discussion

We have explored an example in this section, including a nonlinear singular initial value problem, and a boundary value problem. It is noteworthy that we have developed algorithm in MATLAB code for the proposed LWNN method, and the results have been calculated for these examples.

**Example 1**. We considered Lane-Emden equation as

$$\frac{d^2\xi}{d\eta^2} + \frac{2}{\eta}\frac{d\xi}{d\eta} + \xi^5 = 0,$$

$$\xi(0) = 1, \xi'(0) = 0.$$

And, $\xi(\eta) = \left(1 + \frac{\eta^2}{3}\right)^{-\frac{1}{2}}, \eta \geq 0$ is the exact solution of the example.

Also, the trial solution for LWNN is expressed as
$\xi_t(\eta, p) = 1 + \eta^2 N(\eta, p)$.

In this study, the network was trained using ten equidistant points within the interval [0,1] to compute the results. A comparison between the analytical solutions and the shifted Legendre neural network results is presented in Table 1. These results are also visually compared in Fig. 1, while the error plot is shown in Fig. 2. Table 1 displays the interpolation results at testing points, which were not part of the training data. This testing was conducted to verify if the converged LWNN could accurately produce results when provided with new input points. The comparison shows that the exact (analytical) results closely match the LWNN method results.

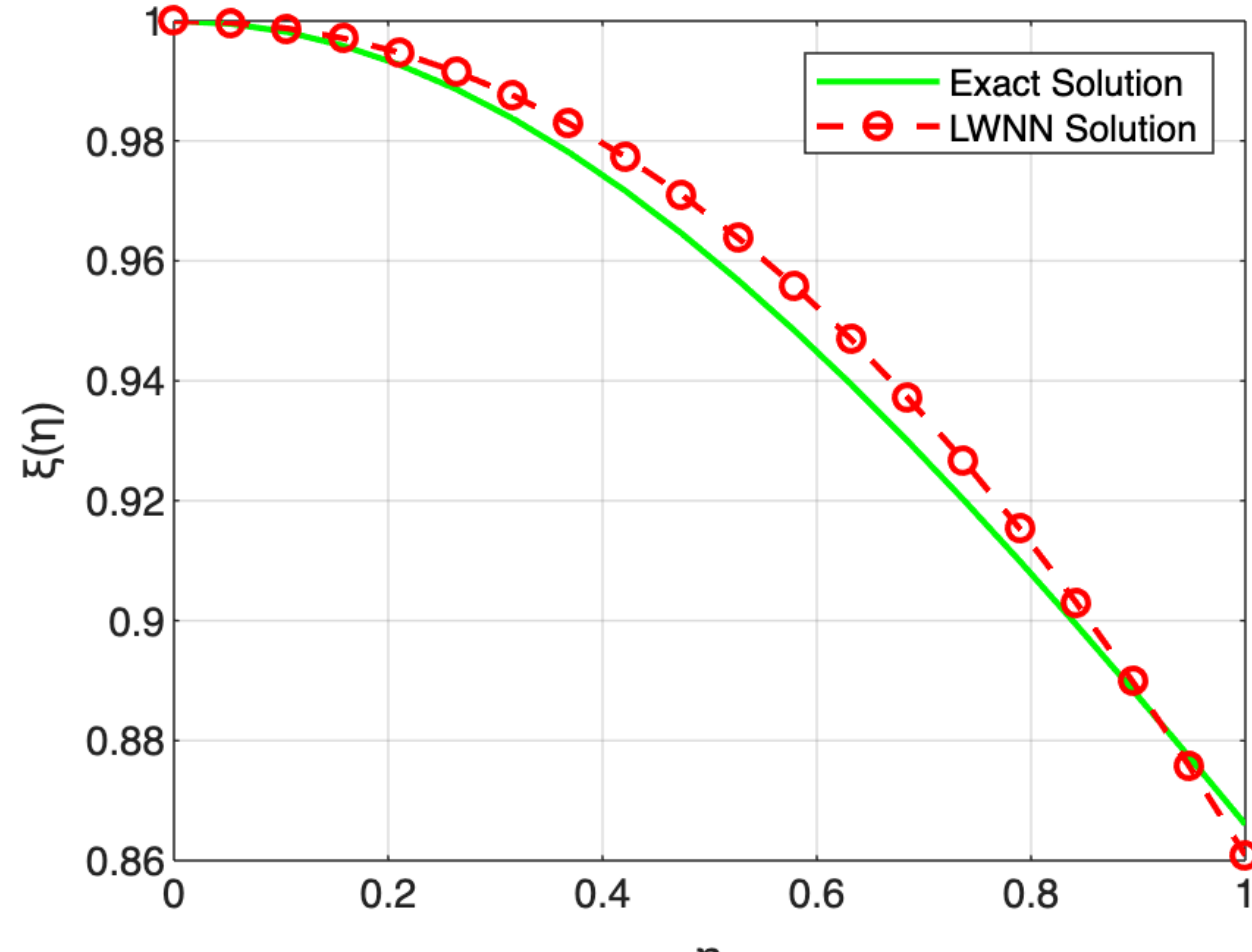

Fig.1. Comparison of exact solution and LWNN

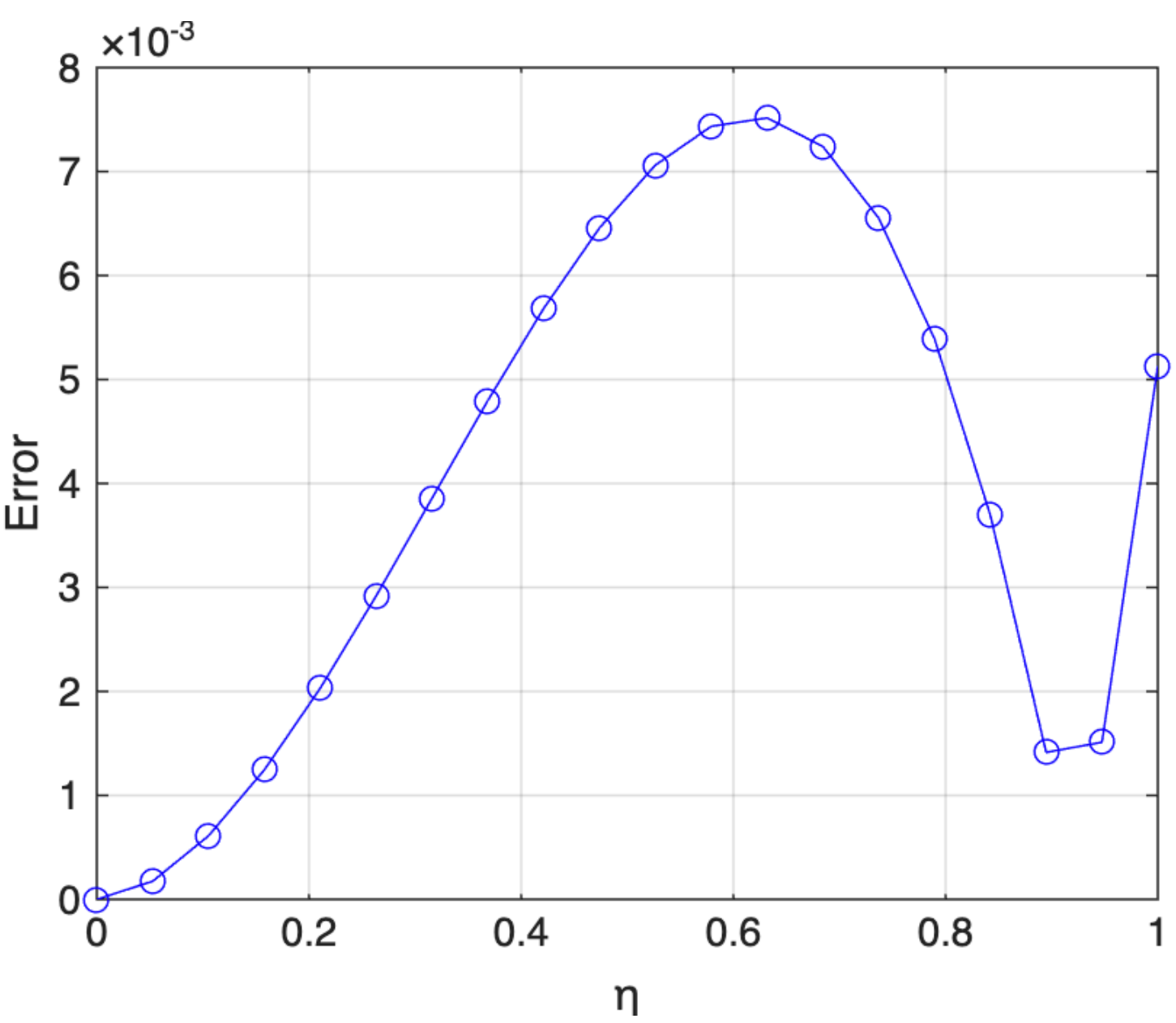

Fig.2. Error between exact solution and LWNN

**Example 2**. We considered Lane-Emden equation as

$$\frac{d^2\xi}{d\eta^2} + \frac{2}{\eta}\frac{d\xi}{d\eta} + \xi = 6 + 12\eta + 2\eta^2 + \eta^3 0 \leq \eta \leq 1,$$

$$\xi(0) = 1, \xi'(0) = 0.$$

And, $\xi(\eta) = \eta^2 + \eta^3, \eta \geq 0$ is the exact solution of the example.

Table 1: Absolute errors history of example 1

| eta | Exact Solution | LWNN Solution | Error |
|---|---|---|---|
| **0** | 1 | 1 | 0 |
| **0.0526** | 0.9995 | 0.9997 | 1.7431E-04 |
| **0.1053** | 0.9982 | 0.9988 | 6.1103E-04 |
| **0.1579** | 0.9959 | 0.9971 | 0.0012 |
| **0.2105** | 0.9927 | 0.9947 | 0.0020 |
| **0.2632** | 0.9887 | 0.9916 | 0.0029 |
| **0.3158** | 0.9838 | 0.9876 | 0.0039 |
| **0.3684** | 0.9781 | 0.9829 | 0.0048 |
| **0.4211** | 0.9717 | 0.9774 | 0.0057 |
| **0.4737** | 0.9646 | 0.9710 | 0.0065 |
| **0.5263** | 0.9568 | 0.9639 | 0.0071 |
| **0.5789** | 0.9484 | 0.9559 | 0.0074 |
| **0.6316** | 0.9395 | 0.9470 | 0.0075 |
| **0.6842** | 0.9301 | 0.9373 | 0.0072 |
| **0.7368** | 0.9202 | 0.9267 | 0.0066 |
| **0.7895** | 0.9099 | 0.9153 | 0.0054 |
| **0.8421** | 0.8993 | 0.9030 | 0.0037 |
| **0.8947** | 0.8885 | 0.8899 | 0.0014 |
| **0.9474** | 0.8773 | 0.8758 | 0.0015 |
| **1** | 0.8660 | 0.8609 | 0.0051 |

Also, the trial solution for LWNN method is expressed as

$$\xi_t(\eta, p) = 1 + \eta^2 N(\eta, p)\,.$$

In this study, the network was trained using ten equidistant points within the interval [0,1] to compute the results. A comparison between the analytical solutions and the shifted Legendre neural network results is presented. These results are visually compared in Fig. 3, while the error plot is shown in Fig. 4. This testing was conducted to verify if the converged LWNN method could accurately produce results when provided with new input points. The comparison shows that the exact (analytical) results closely match the LWNN method results.

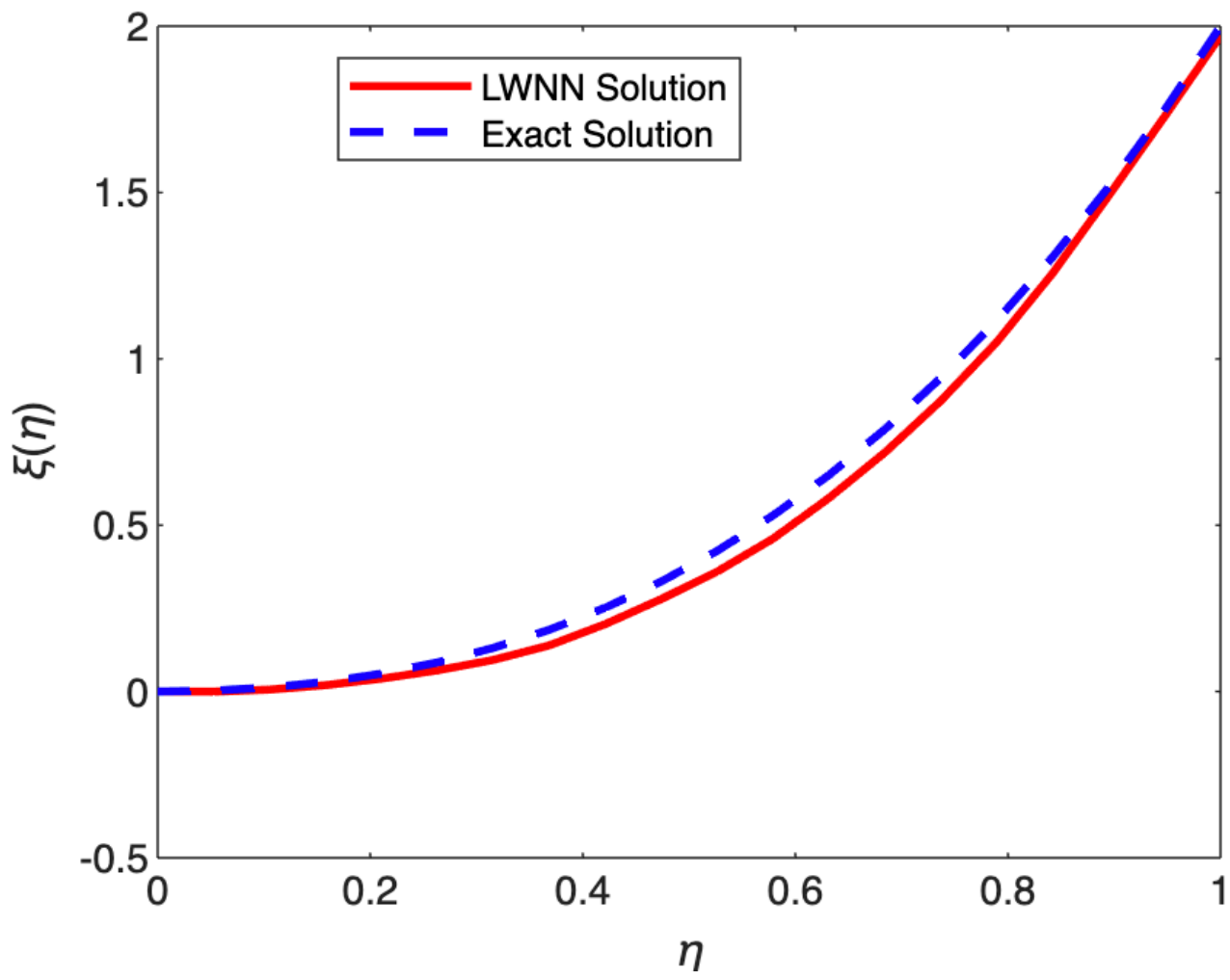

Fig.3. Comparison of exact solution and LWNN solution

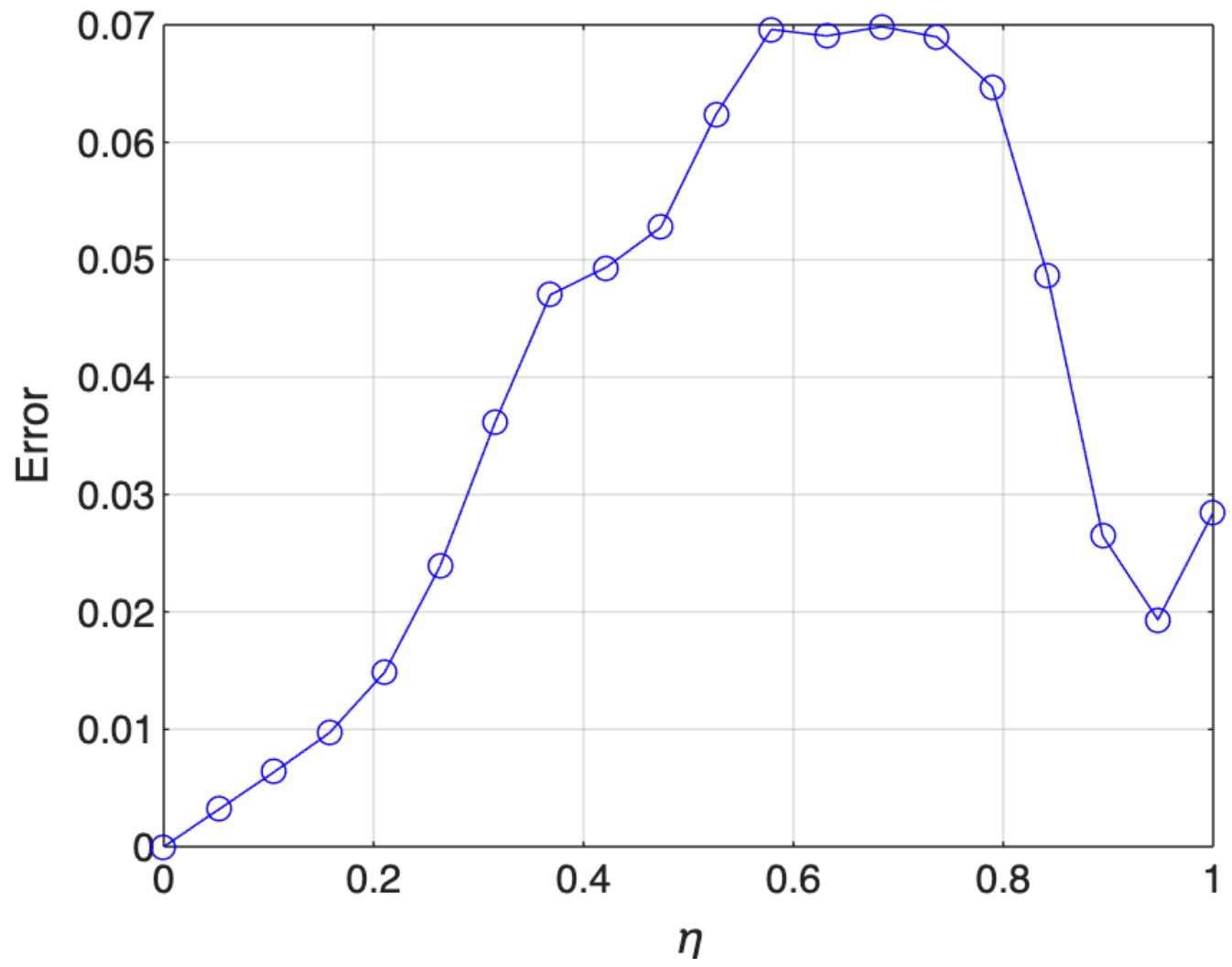

Fig.4. Error between exact solution and LWNN solution

Table 2: Absolute error history of Example 2

| eta | Exact Solution | LWNN Solution | Absolute Error |
|---|---|---|---|
| **0** | 0 | 0 | 0 |
| **0.052632** | 0.0029159 | 0.00027202 | 0.0031879 |
| **0.10526** | 0.012247 | 0.0058807 | 0.006366 |
| **0.15789** | 0.028867 | 0.019172 | 0.0096947 |
| **0.21053** | 0.053652 | 0.038788 | 0.014864 |
| **0.26316** | 0.087476 | 0.06361 | 0.023867 |
| **0.31579** | 0.13121 | 0.095048 | 0.036166 |
| **0.36842** | 0.18574 | 0.13872 | 0.047025 |
| **0.42105** | 0.25193 | 0.20261 | 0.049318 |
| **0.47368** | 0.33066 | 0.2779 | 0.052756 |
| **0.52632** | 0.4228 | 0.36042 | 0.062385 |
| **0.57895** | 0.52923 | 0.45962 | 0.069611 |
| **0.63158** | 0.65082 | 0.58176 | 0.069062 |
| **0.68421** | 0.78845 | 0.7186 | 0.069852 |
| **0.73684** | 0.94299 | 0.87402 | 0.068977 |
| **0.78947** | 1.1153 | 1.0506 | 0.064714 |
| **0.84211** | 1.3063 | 1.2576 | 0.048695 |
| **0.89474** | 1.5168 | 1.4904 | 0.026468 |
| **0.94737** | 1.7478 | 1.7284 | 0.019335 |
| **1** | 2 | 1.9715 | 0.028468 |

## 4. Conclusions

In this paper, we have successfully obtained the application of the Legendre wavelet Neural Network (LWNN) model for solving nonlinear ordinary differential equations second-order of the Lane–Emden type with singularities. By using the Legendre wavelets for functional expansion, we were able to handle the challenges accrued by the singularity in the proposed problem. We have compared the obtain result via LWNN and analytical solutions, and the comparison showed a high degree of accuracy by the low error margins (see Fig 2, Fig 4, Table 1-2). The study also highlighted the advantages of using neural network of a single-layer, where the hidden layer is effectively prevented by expanding the input pattern with Legendre wavelet functions.Moreover, we have implemented the back propagation algorithm for error for updating the network weights proved to be effective in achieving convergence and minimizing the error. The results across various test cases, including nonlinear singular initial value problems, confirmed the robustness and versatility of the LWNN model.

## .Acknowledgements

The authors are very grateful to the editor and reviewer for their constructive comments and suggestions for improvement of the paper